\documentclass[a4paper,12pt]{article}
\usepackage{latexsym}
\usepackage{amssymb}

\usepackage[frame,cmtip,arrow,matrix,line,graph]{xy}

\usepackage{amsmath}
\usepackage{amsthm}
\usepackage{amsfonts}
\usepackage{amscd}

\usepackage{epsfig}

\newtheorem{thm}{Theorem}[section]
\newtheorem{lem}[thm]{Lemma}
\newtheorem{prop}[thm]{Proposition}
\newtheorem{rem}[thm]{Remark}
\newtheorem{ex}[thm]{Example}
\newtheorem{cor}[thm]{Corollary}
\newtheorem{Def}[thm]{Definition}

\newcommand{\C}{\mathcal{C}}
\newcommand{\D}{\mathcal{D}}
\newcommand{\F}{f\mathcal{D}}
\newcommand{\A}{\mathcal{A}}

\newcommand{\B}{\mathcal B}

\newcommand{\Z}{\mathbb Z}

\newcommand{\R}{\mathbb R}

\newcommand{\N}{\mathbb N}

\newcommand{\Top}{\!-\!\textsf{Top}}
\newcommand{\topo}{{\sf Top\;}}

\newcommand{\uh}{\underline{h}}

\newcommand{\al}{\alpha}
\newcommand{\be}{\beta}

\newcommand{\de}{\delta}

\newcommand{\ga}{\gamma}

\newcommand{\te}{\theta}
\newcommand{\s}{\sigma}

\newcommand{\Si}{\Sigma}
\newcommand{\Sn}{\Sigma^n}

\newcommand{\Om}{\Omega}
\newcommand{\On}{\Omega^n}

\newcommand{\lar}{\longleftarrow}
\newcommand{\rar}{\longrightarrow}
\newcommand{\sta}{\stackrel}

\newcommand{\ot}{\otimes}

\newcommand{\x}{\times} 

\newcommand{\vs}{\vspace{5mm}}

\def\z2{\Z_2}

\newcounter{samcounter}
\newenvironment{sam}[3]{
\begin{list}{#2#1{samcounter}#3\stepcounter{samcounter}}{\setcounter{samcounter}{1}}}
{\end{list}}

\begin{document}

\title{Framed disc operads and the equivariant recognition
principle}

\author{Paolo Salvatore \\ {\small{\em Dipartimento di Matematica,
 Universit\`a di Roma ``Tor Vergata''}} \\
{\small{\em Via della Ricerca Scientifica, 
00133 Roma, Italy}}\\
{\small salvator@mat.uniroma2.it}\\
Nathalie Wahl \\
{\small{\em Mathematical Institute, St. Giles' 24-29,}}\\ 
{\small{\em OX1 3LB Oxford, United Kingdom}}\\
{\small wahl@maths.ox.ac.uk}}

\maketitle

\begin{abstract}

The framed $n$-discs operad $f\D_n$ is studied as semidirect
product of $SO(n)$ and the little $n$-discs operad.
Our equivariant recognition principle
says that  a grouplike space acted on by $fD_n$
is equivalent to the $n$-fold loop space on a $SO(n)$-space.
Examples of $fD_2$-spaces are nerves of
ribbon braided monoidal categories.
We compute the rational homology of $fD_n$.
Koszul duality for semidirect product operads
of chain complexes is defined and
applied to compute the double loop space homology
as BV-algebra.
\end{abstract}

\bigskip



\bigskip

\section{Introduction}

The topology of iterated loop spaces was thoroughly investigated in
the seventies \cite{May,Milgram,Segal} .
These spaces have a wealth of homology operations parametrised by the
famous operads of little discs. The notion of operad was introduced
for the first time for this purpose \cite{BV,May}.
Such machinery allows for example
to reconstruct an iterated delooping if one has full knowledge of the
operad action on an iterated loop space.
Moreover any connected space acted on by the little discs is homotopy
equivalent to an iterated loop space. This fact
is the celebrated recognition principle.

Our main objective is to extend this theory by adding
the operations {\em rotating} the discs.
The operad generated by the little $n$-discs $\D_n$
and the rotations in $SO(n)$ 
is called the
framed $n$-discs operad, and
was introduced in \cite{getzler}.
Our equivariant recognition principle (theorem \ref{the})
says that a connected (or grouplike) space
acted on by the framed $n$-discs operad
is weakly homotopic to the $n$-fold loop space on an $SO(n)$-space.

Thus the looping and delooping functors induce a 
categorical equivalence between
$SO(n)$-spaces and spaces acted on by the framed $n$-discs operad, 
under the correct connectivity assumptions.
The main technique consists in presenting the framed
little discs as a semidirect
product of the little discs and the special
orthogonal group,
so that an algebra over the framed $n$-discs is nothing else
than an algebra over the original little $n$-discs
{\em in the category of $SO(n)$-spaces }.

The principle generalises
to any representation $G \to O(n)$ of a topological group.
If $G$ is trivial we recover the original recognition principle
by May.

Is there a way of producing spaces acted on by the
framed discs operad from category theory?
Fiedorowicz showed that the nerve
 of a braided monoidal category is equivalent to an
algebra over the little 2-discs \cite{Fiedo}.
We extend this fact to ribbon braided 
 monoidal categories, which are braided categories equipped with a ``twist'',
and the framed little 2-discs operad (theorem \ref{ribbonthm}).

Is there a procedure to detect framed little 2-discs operad
up to homotopy?
 Fiedorowicz \cite{Fiedo2} has a procedure to recognise
topological operads weakly equivalent to the little 2-discs.
We extend his work to the twisted case giving a procedure to detect
whether an operad is weakly equivalent
to the framed 2-discs (Theorem \ref{fE2}).

\bigskip
We investigate next the homology of spaces acted on by the framed discs operad.
In the category of chain complexes, in analogy with
the topological situation, we define the semidirect product
operad of a Hopf algebra $H$ and an operad acted on by $H$.
An algebra over the semidirect product will be exactly
an algebra over the original operad {\em in the monoidal category
of $H$-modules}.
The homology functor with coefficients in a field commutes with the semidirect product construction.
This approach yields a conceptual proof that the rational
homology of an algebra over the framed 2-discs
is a  Batalin-Vilkovisky algebra \cite{getzler}, 
and computes more generally the rational
homology of the framed $n$-discs operad for any $n$
(Theorem \ref{bvn}).

What is the dual of a semidirect product operad of chain complexes?
We redefine the concept of Koszul duality for semidirect products,
in contrast to \cite{GK}, so that the dual is still naturally
a semidirect product.
This notion induces equivalence of derived categories over dual operads.
Moreover a Koszul operad remains such after taking
semidirect products.
For example the BV-operad is Koszul self-dual up to a shift.
This is based upon the analogous result for Gerstenhaber algebras  \cite{GJ}.

As application we explain how to compute the rational
homology of a double loop
space on a $S^1$-space $X$ as a BV-algebra, starting from a minimal model of
the $S^1$-space together with the derivation induced by the action (Theorem \ref{ratio}).
This extends results in \cite{getzler}, where $X$
is a double suspension,
 and \cite{GJ}, where only the Gerstenhaber algebra structure is
 considered.

We are grateful to Rainer Vogt for organising a stimulating
workshop in Osnabrueck and we thank Tom Leinster, 
Martin Markl, Ulrike
Tillmann and Sasha Voronov
for valuable suggestions.

\section{Semidirect product of topological operads}

We work in the category 
\topo of compactly generated weak Hausdorff topological 
spaces.  
Let $G$ be a topological group. 
The category of left $G$-spaces, denoted $G$-\topo, 
is a symmetric monoidal category by the cartesian product. 
We can thus consider operads in this category, which we
call   {\em $G$-operads}.

Let $\A$ be a $G$-operad. 
 So $\A$ consists 
of a sequence of $G$-spaces $\A(k)$ for $k\in\N$, with $G$-equivariant 
operad structure maps and symmetric group actions.
Note that the unit $1\in\A(1)$ is also preserved by the $G$-action.

We will denote the action of an element $g\in G$ on an element
$a\in\A(k)$ by $ga$.

\begin{Def}\label{defsemi}
Let $\A$ be a $G$-operad. 
Define $\A\rtimes G$,
the \emph{semidirect product of $\A$ and $G$}, to be the following operad 
in \topo: for $k\in\N$, 
\[(\A\rtimes G)(k)=\A(k)\times G^k\]
with $\Si_k$ acting diagonally on the right,
 permuting the components of $G^k$ and acting on $\A(k)$, 
and the map\\
$\ga :
(\A\rtimes G)(k)\times(\A\rtimes G)(n_1)\times\dots\times
(\A\rtimes G)(n_k) \rar (\A\rtimes G)(n_1+\dots+n_k)$\\
 given by

$\ga((a,\underline{g}),(b_1,\uh^1),\dots,(b_k,\uh^k))= 
   (\ga_\A(a, g_1b_1,\dots, g_kb_k),
        g_1 \uh^1,\dots,g_k \uh^k),$\\
where $\uh^i=(h_1^i,\dots,h_{n_i}^i)$ and 
$g_i.\uh^i=(g_i h_1^i,\dots,g_i h_{n_i}^i)$.
The unit in $\A\rtimes G(1)$ is $(1,e)$, formed of the units of 
$\A$ and $G$. 
\end{Def}
The $G$-equivariance of $\ga_A$ is necessary for the associativity of the 
structure map of the semidirect product operad.

\begin{figure}[ht]
\vspace{0.5cm}
\centering
\begin{picture}(0,0)%
\epsfig{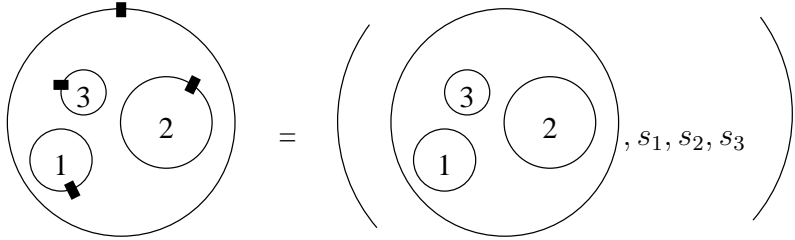}%
\end{picture}%
\setlength{\unitlength}{3947sp}%
\begingroup\makeatletter\ifx\SetFigFont\undefined%
\gdef\SetFigFont#1#2#3#4#5{%
  \reset@font\fontsize{#1}{#2pt}%
  \fontfamily{#3}\fontseries{#4}\fontshape{#5}%
  \selectfont}%
\fi\endgroup%
\begin{picture}(4713,1463)(662,-750)
\put(4546,-151){\makebox(0,0)[lb]{\smash{\SetFigFont{12}{14.4}{\rmdefault}{\mddefault}{\updefault}
$, s_1, s_2, s_3$
}}}
\end{picture}
\caption{element of $f\D_2(3)=\D_2(3)\times(S^1)^3$}
\end{figure}

\begin{ex}\label{framed}
{\rm
The example we have in mind is the framed discs operad $\F_n$. 
Let $\D_n$ be the little n-discs operad of Boardman and Vogt. 
Hence $\D_n(k)$ is the space of embeddings $\coprod_k D^n \to D^n$ of 
$k$ copies of the unit n-disc to itself such that 
the maps are compositions of positive dilations and translations,
and the images are disjoint. 
The framed discs, $f\D_n$, is defined similarly but one allows
rotations for the embeddings.
As spaces, $f\D_n(k)=\D_n(k)\times (SO(n))^k$, the $i$-th
element of $SO(n)$ encoding the rotation of the $i$-th disc. In fact,
$\D_n$ is an $SO(n)$-operad and $f\D_n$ is a semidirect product in the
above sense:
 \[f\D_n = \D_n \rtimes SO(n).\]

The action of $SO(n)$ on $\D_n(k)$ rotates the
little discs around their center. Note that the
whole orthogonal group $O(n)$ acts on $\D_n(k)$ in such a way that 
$\D_n \rtimes O(n)$ is well-defined. We will consider 
semidirect
products of $\D_n$ with any topological group $G$ equipped with a
continuous homomorphism $\phi: G \to O(n)$. We will suppress
$\phi$ from the notation and
denote the resulting semidirect product by $\D_n \rtimes G$.
}
\end{ex}

\begin{prop}\label{spaceACG}
Let $\A$ and $G$ be as in definition \ref{defsemi}. 
A space $X$ is an $(\A\rtimes G)$-algebra if and only if 
$X$ is an $\A$-algebra in the category of $G$-spaces, i.e.
$X$ admits a $G$-action and $\A$-algebra structure maps
$\ \ \te_\A: \A(k)\times X^k \rar X\ $
satisfying \hspace{3mm}
$g(\te_\A(a,x_1,\dots,x_k))=\te_\A( ga, gx_1,\dots, gx_k).$\\
In this case \hspace{3mm}
$\te_{\A\rtimes G}((a,(g_1,\dots,g_k)),x_1,\dots,x_k)=
     \te_\A(a, g_1x_1,\dots, g_kx_k).$
\end{prop}

As an immediate consequence we have

\begin{cor}\label{cor}
Let $X,Y$ be two $(\A\rtimes G)$-algebras. A map $f: X\rar Y$ is a map of
$(\A\rtimes G)$-algebras if and only if it is an $\A$-algebra map 
and a $G$-map.
\end{cor}

We will use the following examples of framed algebras:

\begin{ex}\label{ex}{\rm
Let $Y$ be a pointed $G$-space and 
let $D_n$ denote the monad associated to the
operad $\D_n$ (\cite{May} construction 2.4). $D_n Y$ is the free 
$\D_n$-algebra on the pointed space $Y$.
Let $\On Y$ denote the based $n$-fold 
loop space on $Y$, seen as the space of maps
from the unit $n$-disc $D^n$ to $Y$ sending the boundary to the base point.
The space $\On Y$ carries a natural $\D_n$-algebra structure
\cite{May}.

Let $\phi: G\rar O(n)$ be a continuous group homomorphism. 
The spaces 
$D_n Y$ and $\On Y$ are $\D_n \rtimes G$-algebras, with the action
of $g\in G$

 on $[c;y_1,\dots,y_k]\in D_n Y$, 
where $c\in\D_n(k),\ y_i\in Y$, given by
 \[ g[c;y_1,\dots,y_k]=[\phi(g)c; gy_1,\dots, gy_k],\]

 and on $[y(t)]\in \On Y$, 
where $t\in D^n$ and $[y(t)]$ denotes the n-fold loop $t\mapsto y(t)$, 
given by
\[ g[y(t)]=[gy(\phi(g)^{-1}(t))].\] 
}\end{ex}

\section{Equivariant recognition principle}

Let $\phi: G \rar O(n)$ be as above
and let $X$ be a grouplike $\D_n \rtimes G$-algebra, i.e. the
components of $X$ form a group by the product induced by any element in
$\D_n(2)$. 
As $X$ is a $\D_n \rtimes G$-algebra, it is in particular a
$\D_n$-algebra. 

May introduced a deloop functor $B_n$ from 
$\D_n$-algebras to pointed spaces defined by $B_n X:=B(\Sn,D_n,X)$,
where $B$ is the double bar construction (\cite{May} construction
9.6), and $\Si$ the (reduced) suspension.
May's recognition principle \cite{May,533}
 says that $X$ is weakly equivalent to $\On B_n X$
as $\D_n$-algebra.
May also showed that, conversely, $B_n$ applied to an 
n-fold loop space $\On Y$
produces a space weakly homotopy equivalent to $Y$.

In what follows, we consider the behaviour of $\On$ and $B_n$ with respect to
$G$-actions and provide a recognition principle for algebras over 
$\D_n \rtimes G$. 

\vs

Let $\D_n \rtimes G\Top_{gl}$, $\D_n \rtimes G\Top_0$  and 
$G\Top^*_n$ be the categories of 
grouplike, connected $\D_n \rtimes G$-algebras and $n$-connected pointed
$G$-spaces respectively. Those three categories are closed model 
categories with weak homotopy equivalences as weak equivalences \cite{SV}. 
For a model category $\C$, we will denote Ho($\C$) its associated homotopy 
category, obtained by inverting the weak equivalences.

\vs

For any $G$-space $Y$, we have seen in example \ref{ex} that $\On Y$ has a 
$\D_n \rtimes G$-algebra structure induced by the diagonal action of
$G$. On the other hand, we will define a $G$-action on $B_n X$ for
any $\D_n \rtimes G$-algebra $X$. Hence, $\On$ and $B_n$ will be
functors between the categories of pointed $G$-spaces and of 
$\D_n \rtimes G$-algebras.

\begin{thm}\label{the}
For each continuous homomorphism $\phi: G\rar O(n)$, we have functors
\[\On_\phi=\On: G\Top_{n-1}^* \rar \D_n \rtimes G\Top_{gl}\] 
\[B_n^\phi=B_n:\D_n \rtimes G\Top_{gl} \rar G\Top_{n-1}^*\] 
which induce an equivalence of homotopy categories
\[Ho(G\Top_{n-1}^*) \simeq Ho(\D_n \rtimes G\Top_{gl}).\]
This equivalence restricts to
\[Ho(G\Top_{n}^*) \simeq Ho(\D_n \rtimes G\Top_{0}).\] 
\end{thm}

\begin{proof}
May's recognition principle (\cite{May} Theorem 13.1) is obtained
through the following maps:
\[X\lar B(D_n,D_n,X)\sta{\al}{\rar} B(\On\Sn,D_n,X)\rar 
                                   \On B(\Sn,D_n,X)=\On B_n X,\]
where all maps are $\D_n$-maps between $\D_n$-spaces. 
When $X$ is a $\D_n \rtimes G$-algebra, we want to define $G$-actions on
the spaces involved which induce 
$\D_n \rtimes G$-algebra structures and such that all maps are
$G$-maps. 

The functors
 $D_n$, $\Sn$ and $\On$ restrict to functors in the category of
$G$-spaces, where, for any $G$-space $Y$, we define the action on 
 $D_n Y$, $\Sn Y$ and $\On Y$ diagonally as in example \ref{ex}. 
Hence for any $G$-space $Y$ the $G$-action on $\On\Sn Y$ is given by
\[ g[\s(t),y(t)]=[\phi(g)\s(\phi(g)^{-1}t),gy(\phi(g)^{-1}t)],\]
where $g\in G,\ t,\s(t)\in D^n$ and $y(t)\in Y$.
This produces 
a $\D_n \rtimes G$-algebra structure on $\On\Sn Y$ such that May's map 
$\al:D_n Y \rar \On\Sn Y$ is a $G$-map, and thus a $\D_n \rtimes G$-map.

We extend now these actions on the simplicial spaces $B(D_n,D_n,X)$,
$B(\On\Sn,D_n,X)$ and $\On B(\Sn,D_n,X)$.

Recall that the double bar construction $B(F,C,X)$ is
defined simplicially, for a monad $C$, a left $C$-functor $F$ and
a $C$-algebra $X$ by $B(F,C,X)=|B_*(F,C,X)|$, where $B_p(F,C,X)=FC^pX$, 
with boundary and degeneracy maps using the left
functor, monad and algebra structure maps. The group
$G$ acts then on $B_p(F,C,X)$ through its action on the functors
$F$ and $C$, which comes to ``rotate
everything''.
For example, the action of $g\in G$ on a 1-simplex of $B(\On\Sn,D_n,X)$ is
given by 
$g[\s(t),c(t),x_1(t),\dots,x_k(t)]$\\
$=[\phi(g)\s(\phi(g)^{-1}t), \phi(g)c(\phi(g)^{-1}t), gx_1(\phi(g)^{-1}t),
\dots, gx_k(\phi(g)^{-1}t)]$.

With these actions, all maps above are $G$-maps between 
$\D_n \rtimes G$-spaces and $B_n X$ is
equipped with an explicit $G$-action.

On the other hand, we have a weak homotopy
equivalence  \cite{May, 533}
\[B_n\On Y=B(\Sn,D_n,\On Y)\sta{|\de^*_0|}{\rar} \Sn\On Y \sta{e}{\rar} Y\]
for any $(n-1)$-connected space $Y$. If $Y$ is a $G$-space,
then this composite is a $G$-map with the actions on $B_n\On Y$ and 
$\Si^n\On Y$ defined as above.
\end{proof}

\section{Ribbon braid categories and operads}

We explain in this section how braid groups and ribbon braid groups
are related to the little 2-discs and framed 2-discs operads
respectively. The braid case was done by Fiedorowicz \cite{Fiedo,
Fiedo2}. Details about the ribbon case can be found in 
\cite{transfer,Nthesis}.

We will denote by $\be_k$ the braid group on $k$ strings, the
fundamental group of the configuration space of $k$ unordered
particles in $\R^2$. There is a natural surjection 
$\be_k\twoheadrightarrow\Si_k$, 
sending a braid to the induced permutation of the ends.
The {\em pure braid group}
$P\be_k$ is the kernel of this surjection.

We will denote by $R\be_k $ the
ribbon braid group on $k$ elements,
the fundamental group of the configuration
space of $k$ unordered particles in $\R^2$ with label in $S^1$. 
One can think of an element of $R\be_k$ as a braid on $k$ ribbons, 
where full twists of the ribbons are allowed.
 
The pure ribbon braid group $PR\be_k $ is the kernel of the
surjection $R\be_k \twoheadrightarrow \Sigma_k$.

The groups $R\be_k$ and $PR\be_k$ are isomorphic to 
$\be_k\rtimes\Z^k$  and $P\be_k\times\Z^k$ respectively, where 
$\Z^k$
encodes the number of twists on each ribbon.

We want to characterise operads equivalent to the framed 2-discs.
We consider the following notion of equivalence:

\begin{Def}
An operad map $\A \to \B$ is an {\em equivalence} if
each map $\A(k)\to\B(k)$ is a $\Si_k$-equivariant homotopy equivalence.

An operad $\A$ is a $E_n$-operad (resp. $fE_n$-operad )
if there is a chain of
equivalences connecting $\A$ to $\D_n$
(resp. $f\D_n$).
\end{Def}

Z. Fiedorowicz gave a recognition principle for $E_2$-operads. It
requires the introduction of ``braid operads'', 
which resemble operads except that
the symmetric group actions are replaced by braid group actions
in a natural way. More precisely, a collection of spaces
 $\A=\{\A(k)\}$ is a {\em braid
operad} if $\be_k$ acts on $\A(k)$ for each $k$ and if there are
associative structure maps
$$\gamma: \A(k)\times\A(n_1)\times\dots\times\A(n_k)\to\A(n_1 + \dots +n_k)$$
with two-sided unit $e \in A(1)$, 
satisfying the equivariance conditions

$$\ga(a^\s,b_1,\dots,b_k)=       
 \ga(a,b_{[\s]^{-1}(1)},\dots,b_{[\s]^{-1}(k)})^{\s(n_1,\dots,n_k)}$$ 
and

$$\ga(a,b_1^{\tau_1},\dots,b_k^{\tau_k})=
 \ga(a,b_1,\dots,b_k)^{(\tau_1\oplus\dots\oplus\tau_k)},$$
for all $a\in\A(k), b_i\in\A(n_i), \s\in\be_k, \tau_i\in\be_{n_i},$
where $[\s]$ is the permutation induced by $\s$, the braid  
$\s(n_1,\dots,n_k)$ on $n_1+\dots+n_k$ strings  
is obtained from $\s$ by replacing the ith string by  $n_i$ strings,
and $(\tau_1\oplus\dots\oplus\tau_k)$ 
is the block sum of the braids $\tau_1,\dots,\tau_k$.

A braid operad $\A$ is called a {\em $B_\infty$ operad} if each
$\A(k)$ is contractible and is acted on freely by the braid group
$\beta_k$.

\begin{thm}\label{E2}\cite{Fiedo2}
An operad $\A$ is an $E_2$ operad if and only if its operad structure 
lifts to a
$B_\infty$ operad structure on its universal cover $\tilde{\A}$.
\end{thm}

\noindent
{\em Sketch of the proof.}\hspace{2mm}
Z. Fiedorowicz constructed in \cite{Fiedo} a lift of the operad 
structure of $\D_2$ 
to a $B_\infty$ structure on its universal cover $\tilde{\D_2}$.
The difficulty is that there is no consistent choice of
base-points to lift the map $\ga$. Fiedorowicz uses the inclusion of
the unordered  
little intervals in the little discs $D^u_1$, seeing $D^u_1(k)$ as the 
component of $D_1(k)$ with the intervals ordered in the canonical way, 
from left to right. This provides
contractible subspaces $\D^u_1(k)\subset\D_2(k)$, which will
play the role of base-points for the lifting.

To lift the operad maps for any $E_2$ space, one uses the cofibrant
resolution of $\D_2$, the operad $W\D_2$ constructed by Boardman and
Vogt in \cite{BV} (see also \cite{Vog}).
For any $E_2$ operad $\A$, they show that there is an equivalence of 
operads
$W\D_2\sta{\simeq}{\rar}\A$. We can then use the inclusion
$W\D_1\hookrightarrow W\D_2$ to produce a $B_\infty$ structure 
on $\tilde{\A}$.

For the converse, one first notes that the product of two $B_\infty$
operads is again a $B_\infty$ operad. Now,
using the fact that, as for
$E_\infty$ operads, a braid operad map between two 
$B_\infty$ operads is always an equivalence, the theorem is proved
using the maps
\[\xymatrix{
\tilde{\A} \ar[d] & & 
  \tilde{\A}\times\tilde{\D_2} \ar[ll]_\simeq \ar[rr]^\simeq \ar[d]
          & & \tilde{\D_2}\ar[d]\\
\A & & (\tilde{\A}\times\tilde{\D_2})/P\be_* \ar[ll]_\simeq \ar[rr]^\simeq 
   & & \D_2 .}\] 
\hfill $\square$

\vs

One can define {\em ribbon operads}  and {\em $R_\infty$ operads} 
simply by
replacing the braid groups by ribbon braid groups in the definitions
of braid and $B_\infty$ operads.
As in the braid case, we have that all $R_\infty$ operads are equivalent.

As we have $\D_1\hookrightarrow\D_2\hookrightarrow f\D_2$, 
one can again use the
little intervals to lift the operad maps of $f\D_2$ to its universal
cover $\tilde{f\D_2}$, obtaining then an $R_\infty$ operad
structure, and similarly, for any $fE_2$-operad, using the $W$
construction. 

Adapting the proof of theorem \ref{E2} we obtain:

\begin{thm}\label{fE2}
An operad
$\A$ is an $fE_2$ operad if and only if its operad structure lifts to an
$R_\infty$ operad structure on its universal cover $\tilde{\A}$.
\end{thm}

\vs

Here are our main examples. 
There is a general method to construct categorical operads 
from certain families of groups \cite{T,transfer}. 

The braid groups give rise this
way to an operad $B$, where $B(k)$ is the category with set of objects
$\be_k/P\be_k=\Si_k$ and morphisms 
$Hom_B(\s P\be_k,\tau P\be_k)\cong P\be_k$ 
are by left multiplication 
\[\tau P\be_k \sta{\tau h \s^{-1}}{\lar} \s P\be_k,\]
where $h\in P\be_k$ 
and $\s,\tau \in \be_k$. The operad structure maps 
are defined on object by the usual maps on symmetric groups (from the 
associative operad) and on morphisms by
\[\ga(\s_1\sta{\tau}{\lar}\s_0,\rho_1,\dots,\rho_k)=
      \tau(n_{s_0^{-1}(1)},\dots,n_{s_0^{-1}(k)})
      (\rho_{s_0^{-1}(1)}\oplus\dots\oplus\rho_{s_0^{-1}(k)}),\]
where the right hand side is defined as in the definition of braid 
operads above.
\begin{prop} 
The nerve construction gives an $E_2$-operad $|B|$.  
\end{prop}
\begin{proof}
Let $\tilde{B}(k)$ be the translation category of $\be_k$, having
$\be_k$ as set of objects and 
$Mor_{\tilde{B}(k)}(\s,\tau) = \{\tau\s^{-1} \}$. 
There is then an obvious
``projection'' functor $\tilde{B}\to B$ which induces a covering map
on the nerves. As $|\tilde{B}|$ is contractible, it is in fact the 
universal cover of $|B|$. Also, the operad structure of $|B|$ lifts 
to a natural braid operad structure on $|\tilde{B}|$, which is
a $B_\infty$ operad. We conclude by theorem \ref{E2}.
\end{proof}

Similarly, the ribbon braid groups give rise to a categorical
operad $R$, where the category $R(k)$ has set of objects $\Si_k$ and 
morphisms sets equivalent to $PR\be_k$.
\begin{prop}
The nerve of $R$ yields  an
$fE_2$-operad $|R|$.
\end{prop}
The proof is similar and uses theorem \ref{fE2}. 
Here the universal cover is the nerve of the translation
categories of the groups $R\be_k$.

We note that $|R|$ is the semidirect product $|B|\rtimes B\Z$.

We will describe $B$- and $R$-algebras.

\begin{figure}[ht]
\vspace{0.5cm}
\begin{picture}(0,0)%
\epsfig{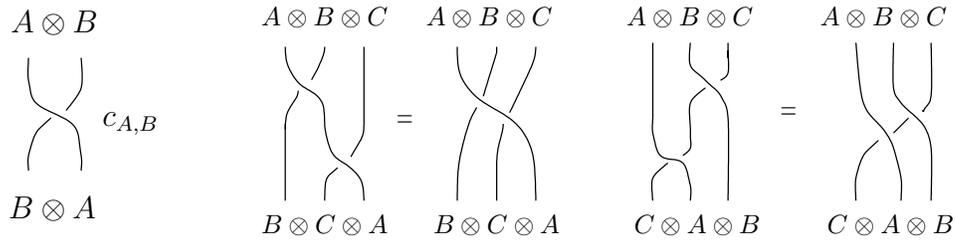}%
\end{picture}%
\setlength{\unitlength}{3947sp}%
\begingroup\makeatletter\ifx\SetFigFont\undefined%
\gdef\SetFigFont#1#2#3#4#5{%
  \reset@font\fontsize{#1}{#2pt}%
  \fontfamily{#3}\fontseries{#4}\fontshape{#5}%
  \selectfont}%
\fi\endgroup%
\begin{picture}(6541,1440)(-34,-736)
\put(541,524){\makebox(0,0)[b]{\smash{\SetFigFont{12}{14.4}{\familydefault}{\mddefault}{\updefault}
$A\otimes B$
}}}
\put(496,-646){\makebox(0,0)[b]{\smash{\SetFigFont{12}{14.4}{\familydefault}{\mddefault}{\updefault}
$B\otimes A$
}}}
\put(811,-61){\makebox(0,0)[lb]{\smash{\SetFigFont{12}{14.4}{\rmdefault}{\mddefault}{\updefault}
$c_{A,B}$
}}}
\put(2206,569){\makebox(0,0)[b]{\smash{\SetFigFont{10}{12.0}{\rmdefault}{\mddefault}{\updefault}
$A\otimes B\otimes C$
}}}
\put(3241,569){\makebox(0,0)[b]{\smash{\SetFigFont{10}{12.0}{\rmdefault}{\mddefault}{\updefault}
$A\otimes B\otimes C$
}}}
\put(2206,-736){\makebox(0,0)[b]{\smash{\SetFigFont{10}{12.0}{\rmdefault}{\mddefault}{\updefault}
$B\otimes C\otimes A$
}}}
\put(3286,-736){\makebox(0,0)[b]{\smash{\SetFigFont{10}{12.0}{\rmdefault}{\mddefault}{\updefault}
$B\otimes C\otimes A$
}}}
\put(4501,569){\makebox(0,0)[b]{\smash{\SetFigFont{10}{12.0}{\rmdefault}{\mddefault}{\updefault}
$A\otimes B\otimes C$
}}}
\put(5716,569){\makebox(0,0)[b]{\smash{\SetFigFont{10}{12.0}{\rmdefault}{\mddefault}{\updefault}
$A\otimes B\otimes C$
}}}
\put(4546,-736){\makebox(0,0)[b]{\smash{\SetFigFont{10}{12.0}{\rmdefault}{\mddefault}{\updefault}
$C\otimes A\otimes B$
}}}
\put(5761,-736){\makebox(0,0)[b]{\smash{\SetFigFont{10}{12.0}{\rmdefault}{\mddefault}{\updefault}
$C\otimes A\otimes B$
}}}
\end{picture}
\caption{braiding and braid relations}\label{braidrel}
\end{figure}

\begin{Def}\label{defbraid}\cite{JS}
A {\em braided monoidal category} is a monoidal category 
$(\A,\otimes)$ equipped with  
a \emph{braiding}, i.e. a natural family of isomorphisms
\[c=c_{A,B} : A\otimes B \rar B\otimes A\]
 satisfying the ``braid relations'' (see figure \ref{braidrel}):\\
$({\rm id}\otimes c_{A,C})\circ a \circ (c_{A,B}\otimes{\rm id})
= a\circ c \circ a : (A\otimes B)\otimes C \to B\otimes(C\otimes A)$,
and\\
$(c_{A,C}\otimes{\rm id})\circ a^{-1} \circ ({\rm id}\otimes c_{B,C})
=a^{-1}\circ c\circ a^{-1}: A\otimes (B\otimes C) \to (C\otimes A)\otimes B$,
where $a$ denotes the associativity isomorphism.

It is called a braided {\em strict} monoidal category if the monoidal
structure is strict, i.e. if the associativity and unit isomorphisms are
given by the identity. 
\end{Def}

Braided monoidal categories arise in the theory of quantum groups and their 
associated link invariants.

\begin{prop}\cite{transfer}
A braided strict monoidal category is exactly
a $B$-algebra.
\end{prop}
If $\A$ is a $B$-algebra, there are functors 
$\te_k:B(k)\times\A^k\to\A$. 
The product in $\A$ is defined on objects by 
$A\otimes B=\te_2(id_{\Si_2},A,B)$ 
and on morphisms by $f\otimes g=\te_2(id_{\be_2},f,g)$, while 
the braiding is given by $c_{A,B}=\te_2(b,id_A,id_B)$, where $b$ is the 
generator of $\be_2$ (see figure \ref{braidrel}). 

The recognition principles of Fiedorowicz and  May
 lead to the following theorem:

\begin{thm}\label{braidthm}\cite{Fiedo}
After group completion
the nerve of a braided monoidal category is
weakly homotopy equivalent to a double loop space.
\end{thm}

\noindent
{\em Sketch of the proof.}\hspace{2mm}
Let $\C$ be a braided monoidal category. Then $\C$ is equivalent to a
braided strict monoidal category $\C'$. Now $|\C'|$ is a
$|B|$-algebra.
Define $X$ to be the double bar construction 
$B(D_2,\tilde{D_2}\times_{P\be}|\tilde{B}|,|\C'|)$. 
 The space $X$ is a $\D_2$-algebra 
weakly equivalent to $|\C'|$. Finally, the recognition
principle tells us that the group completion of $X$ is equivalent to 
a double loop space.
\hfill $\square$

\begin{figure}[ht]
\centering
\begin{picture}(0,0)%
\epsfig{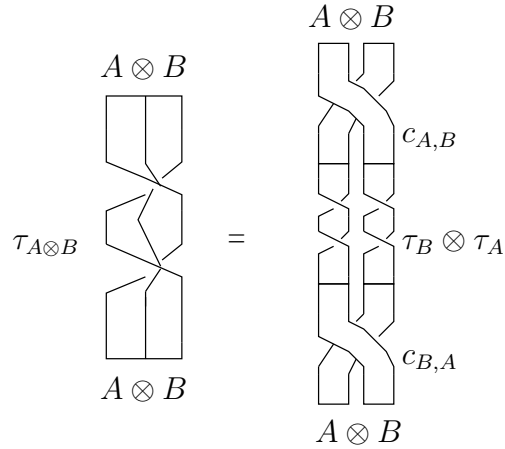}%
\end{picture}%
\setlength{\unitlength}{3947sp}%
\begingroup\makeatletter\ifx\SetFigFont\undefined%
\gdef\SetFigFont#1#2#3#4#5{%
  \reset@font\fontsize{#1}{#2pt}%
  \fontfamily{#3}\fontseries{#4}\fontshape{#5}%
  \selectfont}%
\fi\endgroup%
\begin{picture}(3915,2775)(-1044,-1951)
\put(586,-736){\makebox(0,0)[rb]{\smash{\SetFigFont{12}{14.4}{\rmdefault}{\mddefault}{\updefault}
$\tau_{A\otimes B}$
}}}
\put(991,344){\makebox(0,0)[b]{\smash{\SetFigFont{12}{14.4}{\rmdefault}{\mddefault}{\updefault}
$A\otimes B$
}}}
\put(991,-1681){\makebox(0,0)[b]{\smash{\SetFigFont{12}{14.4}{\rmdefault}{\mddefault}{\updefault}
$A\otimes B$
}}}
\put(2611,-61){\makebox(0,0)[lb]{\smash{\SetFigFont{12}{14.4}{\rmdefault}{\mddefault}{\updefault}
$c_{A,B}$
}}}
\put(2611,-736){\makebox(0,0)[lb]{\smash{\SetFigFont{12}{14.4}{\rmdefault}{\mddefault}{\updefault}
$\tau_B\otimes\tau_A$
}}}
\put(2611,-1456){\makebox(0,0)[lb]{\smash{\SetFigFont{12}{14.4}{\rmdefault}{\mddefault}{\updefault}
$c_{B,A}$
}}}
\put(2296,659){\makebox(0,0)[b]{\smash{\SetFigFont{12}{14.4}{\rmdefault}{\mddefault}{\updefault}
$A\otimes B$
}}}
\put(2341,-1951){\makebox(0,0)[b]{\smash{\SetFigFont{12}{14.4}{\rmdefault}{\mddefault}{\updefault}
$A\otimes B$
}}}
\end{picture}
\caption{compatibility between the twist and the braiding}\label{ribbon}
\end{figure}

\begin{Def}\label{defribbon}\cite{JS}
A {\em ribbon braided (strict) monoidal 
category} $(\A,\otimes,c,\tau)$
is a braided (strict) monoidal  category 
$(\A,\otimes,c)$ equipped with a {\em twist}, i.e. a natural family 
of isomorphisms
\[\tau=\tau_A : A \rar A\]
such that $\tau_1=id_1$, where $1$ is the unit object of $\A$, and 
satisfying the following compatibility with the braiding:
$\tau_{A\otimes B}=c_{B,A}\circ\tau_B\otimes\tau_A\circ c_{A,B}:
A\otimes B\to A\otimes B$
 
(see figure \ref{ribbon}).
\end{Def}

\begin{prop}\cite{transfer}
An $R$-algebra is exactly a ribbon braided strict\\
 monoidal category.
\end{prop}
For an $R$-algebra $\A$, the braided monoidal structure is defined as for
$B$-algebras. 
The twist in $\A$ is defined by $\tau_A=\te_1(t,id_A)$, where $t$ is 
the generator of $R\be_1$.

Consider the monoid 
$\R\x_\Z E\Z \subset \tilde{\F_2}(1)\x_{PR\be_1}|\tilde{R}|(1)$. 
There are monoid maps 
\[S^1\cong (\R\x_\Z *) \ \ \sta{\simeq}{\lar}\  \R\x_\Z E\Z 
                    \ \sta{\simeq}{\rar}\ \  (* \x_\Z E\Z) \cong B\Z.\]
So any $S^1$-space or $B\Z$-space
is canonically an $\R\x_\Z E\Z$-space.
The above maps are restrictions of the operad maps 
$f\D_2 \sta{\simeq}{\lar} \tilde{\F_2}\x_{PR\be}|\tilde{R}| 
\sta{\simeq}{\rar} |R|$ in arity 1. 
Using our recognition principle (theorem \ref{the}) and theorem \ref{fE2}, 
we obtain the following:

\begin{thm}\label{ribbonthm}
The nerve of a ribbon braided monoidal category $\C$, 
after group completion,
 is weakly homotopy equivalent to a 
double loop space $\Om^2 Y$. Moreover, 
the twist on $\C$ induces an $S^1$-action on $Y$ in such a way that 
the equivalence is $\R\x_\Z E\Z$-equivariant.
\end{thm}

\begin{proof} 
Let $\C$ be a ribbon braided monoidal category and let $\C'$ be the 
strictification of $\C$ as a monoidal category.  
The category $\C'$ then inherits a ribbon braided structure 
from the one existing on $\C$.
Its nerve $|\C'|$ is an $|R|$-algebra.
The space $|\C|$ is not necessarily an 
$|R|$-algebra, but it admits a $B\Z$-action induced by the twist on $\C$, 
and the equivalence $|\C|\sta{\simeq}{\rar}|\C'|$ is 
$B\Z$-equivariant.

Now the space 
$X=B(fD_2,\tilde{fD_2}\times_{PR\be}\tilde{|R|},|\C'|)$ is 
weakly homotopy equivalent to $|\C'|$ and is an 
$f\D_2$-algebra. 
The equivalence is obtained through the following diagram of weak equivalences
in $\R\x_\Z E\Z\Top$. 
\[\xymatrix{
B(fD_2,\tilde{fD_2}\times_{PR\be}\tilde{|R|},|\C'|)&  
  B(\tilde{fD_2}\times_{PR\be}\tilde{|R|},
    \tilde{fD_2}\times_{PR\be}\tilde{|R|},|\C'|) \ar[l] \ar[d]  \\
& B(|R|,|R|,|\C'|)\rar |\C'|}\]

The group completion of $X$ is then equivalent to 
a double loop space $\Om^2 Y$, where 
$Y=B(\Si^2,D_2,X)$ and the 
$S^1$-action on $X$ now induces one on $Y$, 
as explained in theorem \ref{the}.
\end{proof}

\section{Semidirect products of algebraic operads}

We work from now on in the category of chain complexes over a field $k$
(possibly with trivial differential). For an element $x$ of a chain
complex, we denote by $|x|$ its degree.
We call operads in this category
differential graded operads, or {\em dg-operads}.

Let $H$ be a graded associative cocommutative Hopf algebra over a field $k$.
The tensor product of two $H$-modules inherits an $H$-structure which is
induced by the coproduct of $H$. As $H$ is cocommutative, the category
of differential graded
$H$-modules, denoted $H$-{\sf Mod}, is a symmetric monoidal category
with product the ordinary tensor product.
Hence it makes sense to consider operads and their
algebras in this category. We call such operads {\em
dg-operads of $H$-modules}.

As in the topological case, we can construct semidirect products for
these operads.

\begin{prop}
Let $P$ be an operad of $H$-modules.
There exists a differential graded operad, the {\em semidirect
product  $P\rtimes H$}, such that algebras over $P$ in the category 
of $H$-modules are exactly $P\rtimes H$-algebras.
\end{prop}

The operad is defined by $(P\rtimes H) (n) = P(n)\ot H^{\ot n}$.
The structure maps are defined similarly to
the topological case, using the comultiplication $c$ of $H$ and 
using interchanging homomorphisms 
with appropriate signs.

\vs

Homology provides a bridge from the
topological to the algebraic setting:
\begin{prop}
Let $G$ be a topological group acting on a topological operad $\A$ .
There is a natural isomorphism 
of operads $H(\A\rtimes G) \cong H(\A)\rtimes H(G).$
\end{prop}

Suppose now that $P$ is a quadratic dg-operad, namely $P$ has 
binary generators  and 3-ary relations  \cite{GK}. 
We will restrict ourselves to the case where $P(1)=k$,
concentrated in dimension 0.
Explicitly $P = F(V)/(R)$, where $F(V)$ is the free operad generated
by a $k[\Si_2]$-module of binary operations $V$ and $(R)$ is the ideal 
generated by a $k[\Si_3]$-submodule $R \subset F(V)(3)$.

\begin{prop}\label{quadop}
Let $H$ be a cocommutative Hopf algebra and\\
 $P=F(V)/(R)$ a 
quadratic operad. 
Then $P$ is an operad of $H$-modules if and only if
\begin{sam}{\roman}{(}{)}
\item $V$ is an $(H,k[\Si_2])$-bimodule;
\item $R \subseteq F(V)(3)$ is an $(H,k[\Si_3])$-sub-bimodule.
\end{sam}
In this case, we will call $P$ a {\em quadratic operad of $H$-modules}.
\end{prop}

\begin{proof}
An element of the free operad on $V$ is described by a tree with vertices labelled
by $V$ \cite{book}.
 
We define the action of $H$ on such  element
by acting on the labels of the vertices, using the comultiplication of $H$.
This is well defined as $H$ is cocommutative. It 
induces an $H$-module structure on $F(V)$ which induces one  
on $P(n)$ for all $n$ by condition (ii). 
The operad structure maps are then $H$-equivariant by construction.
\end{proof}

Let $c:H \to H \ot H$ be the comultiplication.
For $g\in H$
we write informally $(c \ot id)( c(g))= \sum_i g_i'\ot g_i''\ot g'''_i $.

\begin{prop}\label{semiquadratic}
Let $P=F(V)/(R)$ be a quadratic operad of $H$-modules as above.
A chain complex $X$ is an algebra over $P\rtimes H$ if and only if
\begin{sam}{\roman}{(}{)}
\item \label{uno} $X$ is an $H$-module
\item $X$ is a $P$-algebra \label{due}
\item for each $g \in H$, $v \in V$ and $x,y \in X$,  \label{tre}
$$g(v(x,y)) = \sum_i (-1)^{|g_i''||v|+|g_i'''|(|v|+|x|) }
g_i' (v) (g_i''(x), g_i'''(y)).$$
\end{sam}
\end{prop}

\begin{proof}
The $H$-equivariance of the algebra map $\te_2:P(2)\ot X\ot X\rar X$
is given by the commutativity of the following diagram:
\[\xymatrix{
H\ot P(2)\ot X\ot X \ar[rr]^{H \ot \te_2} 
\ar[d]^{{\rm T}\circ (c\ot {\rm id})\circ c} 
                                      & & H\ot X \ar[dd]^\phi \\
H\ot P(2)\ot H\ot X\ot H\ot X \ar[d]^{\psi\ot\phi\ot\phi} & & \\
P(2)\ot X\ot X  \ar[rr]^{\te_2}       & & X, 
} \]
where $\phi$ and $\psi$ give the action of $H$ on $X$ and $P(2)$ 
respectively and $T$ is the interchange.
This diagram translates, for the generators of $P(2)$, into condition (iii)
of the proposition. 
The $H$-equivariance of the structure maps $\te_k$ for $k>2$ is a 
consequence of the fact that $V$ generates $P(k)$, that the operadic 
composition is $H$-equivariant and that the structure maps $\theta$ satisfy
the associativity axiom. 
\end{proof}

\section{Batalin-Vilkovisky algebras}

From now on we work over a field $k$ of characteristic 0.
As first application
we give a conceptual proof of a theorem of Getzler \cite{getzler}.
Recall that a Batalin-Vilkovisky algebra $X$ is
a graded commutative algebra with a linear endomorphism $\Delta:X \to X$
of degree 1 such that $\Delta ^2 =0$ and for each $x,y,z \in X$
the following BV-axiom holds:
\begin{equation} \label{bv}
\begin{split}
\Delta(xyz)= \Delta (xy)z + (-1)^{|x|}x\Delta(yz) +
(-1)^{(|x|+1)|y|}y \Delta(xz) -\Delta(x)yz \\
-(-1)^{|x|}x \Delta (y)z
-(-1)^{|x|+|y|}xy \Delta(z).
\end{split}
\end{equation}

\begin{thm}\label{bvex}\cite{getzler}
Let $H(fD_2)=H(D_2) \rtimes H(SO(2))$ be the homology of the
framed little 2-discs operad. An $H(fD_2)$-algebra is exactly 
a Batalin-Vilkovisky algebra. 
\end{thm}

\begin{proof}
Let $X$ be an algebra over $H(fD_2)$. 
By proposition \ref{semiquadratic} (condition (i)), 
$X$ is an $H(SO(2))$-module.
As an algebra, $H(SO(2))=k[\Delta]/\Delta^2$, where $\Delta \in H_1(SO(2))$
is the fundamental class. 
This provides $X$ with 
an operator $\Delta$ of degree 1 satisfying $\Delta^2=0$.
Condition (ii) of proposition \ref{semiquadratic}
tells us that $X$ is an algebra over $H(D_2)$.
The operad $H(D_2)$, called the Gerstenhaber operad,
was identified by F. Cohen \cite{cohen}.
This operad is
quadratic, generated by the
operations $* \in H_0(D_2(2))$ and $b \in H_1(D_2(2))$, corresponding
to the class of a point
and the fundamental class under the $SO(2)$-equivariant
homotopy equivalence $D_2(2) \simeq S^1$.
The class $*$ induces 
a graded commutative product on $X$, while 
$b$ induces a Lie bracket of degree 1, i.e. a 
Lie algebra structure 
on $\Si X$, the suspension of $X$, where $(\Si X)_i=X_{i-1}$.
The bracket is defined on $X$ by 
$[x,y]=(-1)^{|x|}b(x,y)$.  Cohen proved that 
the product and the bracket satisfy the following Poisson relation:
\begin{equation} \label{poisson}
[x,y*z] = [x,y]*z +(-1)^{|y|(|x|+1)} y*[x,z]  \, .
\end{equation}

In order to unravel condition (iii) of proposition \ref{semiquadratic}, 
 we must understand the effect in
homology of the $SO(2)$-action on $D_2(2) \simeq S^1$. 
Clearly  $\Delta(*)=b$ because
the rotation of the generator in degree 0 gives
precisely the fundamental 1-cycle. Moreover $\Delta(b)=0$ for dimensional 
reasons.
As $\Delta$ is primitive, condition (iii)  
applied respectively to 
$g=\Delta, \, v=*$ and $g=\Delta,\, v=b$ provides the following relations:
\begin{equation} \label{quattro}
 \Delta(x*y) =  \Delta(*)(x,y) + \Delta(x)*y + (-1)^{|x|} x* \Delta(y)
 \, ;
\end{equation}
\begin{equation} \label{cinque}
 \Delta(b(x,y)) = \Delta(b)(x,y) - b(\Delta(x),y)
 +(-1)^{|x|+1}b(x,\Delta(y))\, .
\end{equation}

As $\Delta(*)=b$,
equation \ref{quattro} expresses
the bracket in terms of the product and $\Delta$ :
\begin{equation} \label{sei}
[x,y]=(-1)^{|x|}\Delta(x*y)-(-1)^{|x|}\Delta(x)*y - x * \Delta(y) \, .
\end{equation}

If we substitute this expression into the Poisson relation \ref{poisson}
we get exactly the BV-axiom \ref{bv}. We can re-wright 
equation \ref{cinque} as
\begin{equation}
\Delta[x,y]=[\Delta(x),y]+(-1)^{|x|+1}[x,\Delta(y)]
\end{equation}
which  says that $\Delta$ is a
derivation with respect to the bracket. 
To conclude we must show that \ref{cinque} and the Lie 
algebra axioms follow from the BV-axiom. 
This is shown in Proposition 1.2 of \cite{getzler}.
\end{proof}

\begin{figure}[ht]
\begin{center}
\begin{picture}(0,0)%
\epsfig{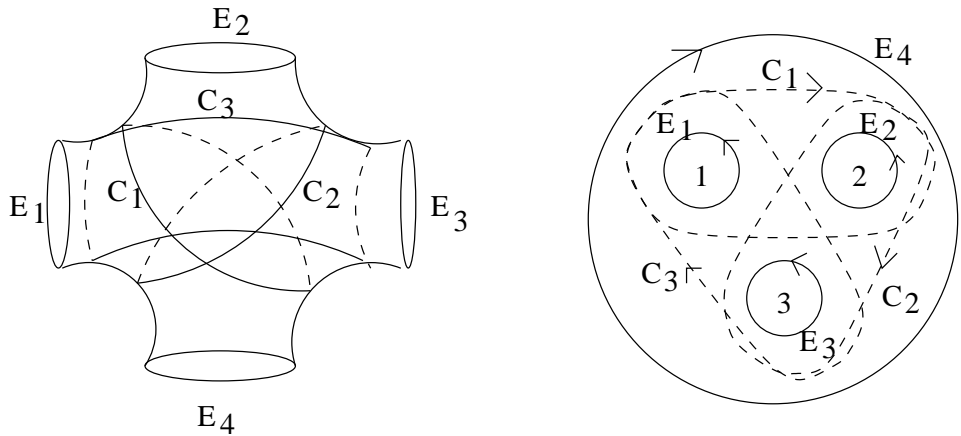}%
\end{picture}%
\setlength{\unitlength}{3947sp}%
\begingroup\makeatletter\ifx\SetFigFont\undefined%
\gdef\SetFigFont#1#2#3#4#5{%
  \reset@font\fontsize{#1}{#2pt}%
  \fontfamily{#3}\fontseries{#4}\fontshape{#5}%
  \selectfont}%
\fi\endgroup%
\begin{picture}(5680,2565)(136,-2266)
\end{picture}
\caption{Lantern relation}\label{lantern}
\end{center}
\end{figure}

\begin{rem}{\rm
The lantern relation, introduced by Johnson for its relevance to the
mapping class group of  surfaces \cite{J} is defined by the following
equation: $T_{E_4}=T_{E_1}T_{E_2}T_{E_3}T_{C_1}T_{C_2}T_{C_3}$, where
$T_C$ denotes the Dehn twist along the curve $C$. See figure
\ref{lantern} for the relevant curves on a sphere with four holes, or
equivalently on a disc with three holes. 
The mapping class group of a sphere with 4 ordered holes, 
is the group of path
components of orientation preserving diffeomorphisms which fix
the boundary pointwise.
The group is isomorphic to the pure ribbon braid group $PR\be_3$.
The lantern relation is thus a
relation in  $PR\be_3$ and gives rise to a
relation in $H_1(f\D_2(3))$ which is the abelianisation of $PR\be_3$. 
It was noted by Tillmann that, with this
interpretation, one gets precisely the BV-axiom \ref{bv}. Indeed, 
up to signs, the curve $E_1$ represents
the operation $(x,y,z) \mapsto 
\Delta x*y*z$. Moreover 
$E_2$ corresponds to $x*\Delta y*z$, $E_3$ to $x*y*\Delta z$,
$E_4$ to $\Delta(x*y*z)$, $C_1$ to $\Delta(x*y)*z$, $C_2$ to $x*\Delta(y*z)$ 
and $C_3$ to $y*\Delta(x*z)$.

This geometric interpretation shows that any
$H(f\D_2)$-algebra is a BV-algebra.}
\end{rem}

\bigskip

We will use alternatively the notations $e_n$ and $H(D_n)$ for the 
homology of the little $n$-discs operad, by which we mean
$$\left\{ \begin{array}{ll}
e_n(k):=H(\D_n(k)) & k\ge 1\\
e_n(0)=0.           &  
\end{array}\right.$$ 

Algebras over the operad $e_n$,
$n\ge2$,  are called
{\em $n$-algebras}. By assumption they have no units. 
F. Cohen's study of $H(\D_n)$ in \cite{533}
implies that 
an $n$-algebra $X$ is a differential graded
commutative algebra with a Lie bracket of degree $n-1$, i.e. 
\begin{sam}{\arabic}{($L$}{)}
\item $[x,y]+(-1)^{(|x|+n-1)(|y|+n-1)}[y,x]=0,$
\item $\partial[x,y]=[\partial x,y]+(-1)^{|x|+n-1}[x,\partial y]$, 
\item $[x,[y,z]]=[[x,y],z]+(-1)^{(|x|+n-1)(|y|+n-1)}[y,[x,z]]$,
\end{sam}
satisfying the Poisson relation 
\begin{sam}{\arabic}{($P$}{)}
\item $[x,y*z] = [x,y]*z +(-1)^{|y|(|x|+n-1)} y*[x,z]  \, .$
\end{sam}
Gerstenhaber algebras correspond to the case $n=2$.

Note that
$\D_n(2)$ is $SO(n)$-equivariantly homotopic to $S^{n-1}$.
In an 
$n$-algebra,  
the product comes from the generating class $* \in H_0(\D_n(2)) \cong k$ 
and the bracket from the fundamental class $b \in 
H_{n-1}(\D_n(2)) \cong k$, if we define $[x,y]=(-1)^{(n-1)|x|}b(x,y)$.
The operad $e_n$ is quadratic \cite{GJ}.

In order to determine the homology operad $H(f\D_n)$, we need to
know the Hopf algebra structure of $H(SO(n))$ and
the effect in homology of the action of $SO(n)$ on $\D_n(2)$. 
For dimensional reasons, one always has $\de(b) =0$ for
each $\de \in \tilde{H}(SO(n))$.
On the other hand $\de(*)=\pi_*(\delta)$, where $\pi_*$ is induced
in homology by the
evaluation map $\pi: SO(n) \to S^{n-1}$, via
$\D_n(2) \simeq S^{n-1}$.

Let us give two further examples:

\begin{ex}\label{fD3-4}
(i) An $H(f\D_3)$-algebra 
is a $3$-algebra together with an endomorphism $\delta$ of degree 3
such that $\delta^2 =0$ and $\delta$ is a derivation  both with respect 
to the product and the bracket.

(ii) An $H(f\D_4)$-algebra is a commutative dg-algebra
together with two linear endomorphisms $\alpha, \beta$
of degree $3$, such that $\alpha^2 =0, \, \beta^2=0, \, 
\alpha\beta = - \beta \alpha$, 
$\alpha$ and the product satisfy the BV-axiom \ref{bv}, and
$\beta$ is a derivation with respect to the product. 
\end{ex}

\begin{proof}
(i) $H(SO(3))=\wedge(\de)$ is the free exterior algebra generated by the
fundamental class
$\delta \in H_3(SO(3))$, and
$\pi_*(\delta)=0$ by dimension.   
By proposition \ref{semiquadratic}, $X$ is an $H(f\D_3)$-algebra if
and only if $X$ is a 3-algebra admitting 
an $H(SO(3))$-module structure, i.e. an operator $\de$
of degree 3 with $\de^2=0$, such that the following
relations hold:
\begin{equation}\label{de*}
\de(x*y)=\de x*y+(-1)^{|x|}x*\de y 
\end{equation}
\begin{equation}\label{deb}
\de[x,y]=[\de x,y]+(-1)^{|x|}[x,\de y] \, , 
\end{equation}
where those equations are obtained by setting $g=\de, \, v=*$ and 
$g=\de, \, v=b$ in turn in condition (iii) of proposition 
\ref{semiquadratic}. In this case the operator $b$ is
equal to the bracket and lies in degree 2 and $\de(b)=0$.

(ii) The evaluation fibration $SO(3) \to SO(4) \to S^3$ splits as a product.
So $H(SO(4))=\wedge(\al,\be)$, with both generators in degree 3. 
A class  $\alpha$ comes from the basis, so 
$\pi_*(\alpha)=b$, and
another class $\beta$ comes from the fibre, so $\pi_*(\beta)=0$.

As in the previous case, we know that an $H(f\D_4)$-algebra $X$ is a
4-algebra with two operators $\al$ and $\be$ both in degree 3,
satisfying $\al^2=0=\be^2$ and $\al\be=-\be\al$ and  relations 
obtained by settting $(g,v)=(\al,*)$,  $(g,v)=(\al,b)$,
$(g,v)=(\be,*)$, and 
$(g,v)=(\be,b)$ in turn in condition (iii) of proposition \ref{semiquadratic}. 
Using the identification $b(x,y)=(-1)^{|x|}[x,y]$, this gives the 
following equations:
\begin{equation}\label{al*}
\al(x*y)=(-1)^{|x|}[x,y]+\al x*y+(-1)^{|x|}x*\al y 
\end{equation}
\begin{equation}\label{alb}
\al[x,y]=[\al x,y]+(-1)^{|x|+1}[x,\al y] 
\end{equation}
\begin{equation}\label{be*}
\be(x*y)=\be x*y+(-1)^{|x|}x*\be y 
\end{equation}
\begin{equation}\label{beb}
\be[x,y]=[\be x,y]+(-1)^{|x|+1}[x,\be y] \, . 
\end{equation}
Note that equations \ref{al*} and \ref{alb} correspond precisely to
the equations we had for $\Delta$ and the bracket in theorem
\ref{bvex}. So, by the same calculations, we know that $\al$ and the
product form a Batalin-Vilkovisky algebra of higher degree, i.e. $\al$
and $*$ satisfy equation \ref{bv} but the operator $\al$ has now 
 
degree 3. There is an additional operator $\be$ of degree 3. Equation
\ref{be*} says that $\be$ is a derivation with respect to the
product. Using equation \ref{al*}, one can rewrite equation \ref{beb} in
terms of $\al, \be$ and the product. This shows that equation
\ref{beb} is redundant.
\end{proof}

We need a lemma in order to state the general case.
\begin{lem}\label{SOn}
For $n\ge 1$, over a field of characteristic 0, the Hopf algebra  
$H(SO(2n))=\bigwedge(\be_1,\dots,\be_{n-1},\al_{2n-1})$ 
is the free exterior algebra on primitive generators
$\be_i \in H_{4i-1}(SO(2n))$
and $\al_{2n-1} \in H_{2n-1}(SO(2n))$. Moreover, $\pi_*(\be_i)=0$ 
for all $i$
and $\pi_*(\al_{2n-1})=b \in H_{2n-1}(S^{2n-1})$ is the fundamental class. 

The Hopf algebra $H(SO(2n+1))=\bigwedge(\be_1,\dots,\be_n)$ is the free 
exterior algebra on
primitive generators $\be_i \in H_{4i-1}(SO(2n+1))$, 
and $\pi_*(\be_i)=0$ for all $i$.
\end{lem}

\begin{proof}
The homology Serre spectral sequence of the principal fibration\\
$SO(n) \to SO(n+1) \to S^{n}$
collapses at the $E_2$ term if $n$ is odd;
if $n$ is even then there is a non-trivial differential
$d(b)= \al_{n-1}$ \cite{Mimura}.
\end{proof}

If a Hopf algebra $H$ acts trivially,
via the counit, on an operad $P$, we call 
the semidirect product just the {\em direct product} and denote it  by
$P \times H$.
Note that a $P\times H$-algebra is an $H$-module $X$ with a
$P$-algebra structure satisfying an $H$-equivariance condition which
is trivial only if $H$ acts trivially on $X$. In particular, any
$P$-algebra is a $P\times H$-algebra with the trivial $H$-module
structure.

Let us denote by $BV_n$, for $n$ even, the
Batalin-Vilkovisky operad with the operator  $\Delta$
in degree $n-1$. Hence a $BV_n$-algebra is a differential graded
commutative algebra with an operator $\Delta$ of degree $n-1$ such that 
$\Delta^2=0$ and the BV-equation (\ref{bv}) holds.

Note that there is no non-trivial $\Si_2$-equivariant map from 
$H_0(\D_{2n+1}(2))$ to $H_{2n}(\D_{2n+1}(2))$. So $(2n+1)$-algebras 
do not give rise to ``odd'' $BV$-structures like in the even case.

\begin{thm} \label{bvn}
For $n\ge 1$ there are isomorphisms of operads
\[H(f\D_{2n+1})\cong H(\D_{2n+1})\times H(SO(2n+1))\]
and 
\[H(f\D_{2n})\cong BV_{2n}\times H(SO(2n-1)).\]
\end{thm}

Hence an $H(f\D_{2n+1})$-algebra is a $(2n+1)$-algebra together with
endomorphisms
$\be_i$ of degree $4i-1$ for $i=1,\dots,n$ such that
$\be_i^2=0$,  $\; \be_i\be_j=-\be_j\be_i$ for each $i,j$, and each $\be_i$ is
a $(2n+1)$-algebra derivation, i.e. 
a derivation both with respect to the product and the bracket.

On the other hand, an $H(f\D_{2n})$-algebra is a $BV_{2n}$-algebra together
with endomorphisms $\be_i$ of degree $4i-1$ for $i=1,\dots,n-1$
squaring to 0 and anti-commuting as in
the odd case, 
which moreover anti-commute with the $BV$ operator $\Delta$ and are 
derivations with respect to the product.

The proof is similar to the proof of example \ref{fD3-4}.

We already saw that iterated loop spaces are algebras over the framed
discs operad. We deduce the following example:
\begin{ex}
The homology of an $n$-fold loop space on a pointed $SO(n)$-space 
is an algebra over  $H(f\D_n)$.
\end{ex}

Another interesting class of algebras over $H(fD_n)$ is given by 
the space $\Lambda^n(X)$ of unbased maps from
$S^n$ to a space $X$.
Chas and Sullivan showed that the homology of a free loop space
$\Lambda M$ on
an oriented  manifold $M$ is a Batalin-Vilkovisky algebra. 
Sullivan and Voronov generalised it to higher dimension and have a 
geometrical proof involving the so-called cacti operad.

\begin{ex}\cite{SS} 
Let $M$ be a $d$-dimensional oriented manifold.
Then the $d$-fold desuspended homology
$\Si^{-d} H(\Lambda^n M)$
of the unbased mapping space
from the $n$-sphere into $M$
is an algebra over $H(fD_{n+1})$. 
\end{ex}

\section{Koszul duality for semidirect products}

Recall that we work in the category of differential graded vector 
spaces, also called chain complexes, 
over a field $k$ of characteristic 0. 
Let $P$ be a quadratic operad of $H$-modules. 
We assume that 
$P(0)=0$ and $P(1)=k$ is concentrated in degree 0. 

Recall that if $P=F(V)/(R)$  
is the quadratic operad generated by $V$ with
relations $R$, then
its quadratic dual, as defined in \cite{GK} and
\cite{book}
 is given by
$P^!:=F(\check{V})/(R^\bot)$, 
where $\check{V}=V^*\ot sgn_2$, $V^*$ is the linear dual,
$sgn_2$ is the sign representation of the symmetric group, 
and $R^\bot$ is the annihilator of $R$ in $F(\check{V})(3)$.  

The dual $(P\rtimes H)^!$ of $P\rtimes H$ in this sense is not
naturally a semidirect product operad. Thinking of $P\rtimes H$ as the 
operad $P$ in $H$-{\sf Mod}, we consider instead the following duality.

\begin{Def}
The {\em dual of a semidirect product} $P\rtimes H$ is the operad \\
$(P\rtimes H)^\dagger:=P^!\rtimes H^{op}$, 
where $P^!$ is the quadratic dual of $P$ and
$H^{op}$ denotes $H$ with the opposite multiplication.
\end{Def}

This makes sense because
$P^!$ is an operad of
$H^{op}$-modules by proposition \ref{quadop}. 

\vs

The suspension of an operad $P$
is the operad $\Lambda P$
defined by
$\Lambda P(n)=\Si^{n-1}(P(n))\ot sgn_n$ \cite{GJ}.
A chain complex $A$ is
a $\Lambda P$-algebra if and only if the suspension
$\Si A$ is a $P$-algebra.

In the following example, we show that the operad $BV$, as a 
semidirect product, is self-dual up to suspension.

\begin{ex}
\begin{enumerate} 
\item $BV_{2n}^\dagger:=H(\D_{2n})^!\rtimes k<\al_{2n-1}>
                       =\Lambda^{1-2n}BV_{2n}$;
\item $H(f\D_{2n})^\dagger:=H(\D_{2n})^!\rtimes H(SO(2n))^{op}=
                \Lambda^{1-2n}BV_{2n}\times H(SO(2n-1))$;
\item $H(f\D_{2n+1})^\dagger= \Lambda^{-2n}H(\D_{2n+1}) \times H(SO(2n+1))$.
\end{enumerate}
\end{ex}

\begin{proof}
We give a proof of (1).  It is known that the quadratic
dual of the operad $e_{2n}=H(\D_{2n})$
 is its own $(2n-1)$-fold desuspension $\Lambda^{1-2n} e_{2n}$,   
with product $p'=b^*$ dual to original bracket $b$ and bracket
$b'=p^*$ dual to the original product $p$ \cite{GJ, Markl}.
If $\al_{2n-1}(p)=b$, then $\al_{2n-1}(b^*)=p^*$. So the class 
$\al_{2n-1}$ gives an operator $\Delta$ of degree $2n-1$ with 
 $\Delta(p')=b'$ and $\Delta(b')=0$, which thus induces a $BV$ structure 
as in theorem \ref{bvex} but this time with product in degree $1-2n$.
\end{proof}

From now on, we assume that $P$ is a quadratic operad of 
$H$-modules such that 
 $P(n)$ is a finite dimensional $k$-vector spaces for each $n$.
The operad $P$ is
thus {\em admissible} in the sense of Ginzburg and Kapranov.

Let $P-\textsf{Alg}_{\geq n}$
($P-\textsf{Alg}_{\leq n}\,$ ) denote the category
of $P$-algebras of finite type concentrated in degree $\geq n$ ($\leq n$).
Getzler and Jones (\cite{GJ}, see also \cite{muriel}) 
constructed contravariant adjoint functors
$$C_P: P-\textsf{Alg}_{\geq 1}   
             \rightleftarrows P^!-\textsf{Alg}_{\leq -2} :T_{P} $$
such that                           
$C_P$ and $T_{P}$ preserve quasi-isomorphisms. Moreover,
 if $P$ is Koszul, the unit and counit
of the adjunction are quasi-isomorphisms.

\

If $\C$ is a category of chain complexes, we denote by $Ho(\C)$ the category
obtained from $\C$ by inverting all quasi-isomorphisms.
We thus have the following equivalences of categories:
\[Ho(P-\textsf{Alg}_{\geq 1}) \simeq
           Ho(P^!-\textsf{Alg}_{\leq -2})   \]

We want to see that $C_P$ and $T_P$ restrict to functors
between categories 
of  $P\rtimes H$-algebras 
and  $P^!\rtimes H^{op}$-algebras.

Let $A$ be a $P$-algebra in $H$-{\sf Mod}. 
The complex $C_P(A)$ is defined as the 
free $P^!$-algebra on the
dual of the suspension of $A$
with differential $d_1+d_2$, where $d_1$ is induced by the 
differential of $A$ and $d_2$ by the $P$-algebra structure of $A$. 
This object has an $H$-module structure induced by the action of $H$ on $A$ and on 
$P$. 
The $H$-action commutes with $d_1$ as $A$ is an $H$-module and it 
commutes with $d_2$ as the $P$-algebra structure maps of $A$ are $H$-equivariant. 
The $P^!$-algebra structure map of
$C_P(A)$
 is $H$-equivariant because $P$ is an operad of
$H$-modules.

Similarly, for a $P^!$-algebra $X$,  
$T_P(X)$ is the free $P$-algebra on the dual of the suspension
of $X$ with twisted differential,
and is a $P$-algebra in $H-{\sf Mod}$.
Moreover the natural transformations $T_P C_{P} (A) \to A$ and
$C_{P} T_P (X)  \to  X$ are respectively
$H$- and $H^{op}$-equivariant.

We thus have the following theorem:

\begin{thm}
Let $P$ be a quadratic operad of $H$-modules. Then there
is an equivalence of
categories
\[Ho(P\rtimes H-\textsf{Alg}_{\geq 1}) \simeq 
           Ho(P^!\rtimes H^{op}-\textsf{Alg}_{\leq -2}). \]
\end{thm}

\begin{rem}{\rm 
The theorem can be extended to objects not of finite type,
by considering coalgebras over an operad.
}\end{rem}

\begin{Def}
Let $A$ be a $P\rtimes H$-algebra. The {\em operadic homology of A} is 
the $H(P^!)\rtimes H^{op}$-algebra $H(C_P(A))$. 
\end{Def}
So the homology of $A$ over $P\rtimes H$ is the 
homology of $A$ over $P$ \cite{book}
with an additional 
$H^{op}$-module structure.

We apply this machinery
in order to compute explicitly the $BV$-algebra
 structure of the homology of a double loop space.
The Gerstenhaber algebra structure is computed
in 6.1 of \cite{GJ} for the double loop space on a manifold 
via operadic homology.

Let $X$ be a $2$-connected pointed  CW-space of finite type acted on by $S^1$.
Let $M(X)$ be the minimal model of $X$ \cite{FeHaTh}, a non-negatively
graded unital commutative algebra with a differential of degree 1. 
The $S^1$-action induces a map $f: M(X) \to M(X) \ot \Lambda(e_1)$ of
type $f(x)=x + \Delta(x) \ot e_1$, thus $\Delta$ is a derivation of degree
-1. Let $m(X)=M(X)/ <1>$ be the quotient by the unit.
By changing signs we regard $m(X)$ as a BV-algebra concentrated in degree
$\leq -3$, with trivial Lie bracket. Clearly the suspension
$\Si m(X)$ is a $BV^\dagger$-algebra concentrated in negative degree.

\begin{thm} \label{ratio}
The operadic homology of the $BV^\dagger$-algebra $\Si m(X)$  
is isomorphic to $H(\Omega^2(X))$ as $BV$-algebra.
\end{thm}

\begin{proof} 

If we substitute the expression $BV$ by $G$, then the result
follows from 6.9 of \cite{GJ} by the same proof as 6.1 of \cite{GJ},
with the minimal model replacing the deRham complex.

Since $H_*(\Omega^2 X)$ is the free commutative algebra on 
$\pi_*(\Omega^2 X)$, by the Milnor-Moore theorem, 
the homology $BV$-algebra structure is uniquely determined
by $\Delta$ on spherical classes.

We identify geometrically the $BV$-operator
on spherical classes.
If $Y$ is a pointed $S^1$-space and $x \in \pi_n(Y)$, consider the map
$S^n \wedge S^1_+ = (S^n \times S^1)/(* \times S^1) \to Y$
induced by the action.
Since $S^n \wedge S^1_+ \simeq S^n \vee S^{n+1}$ let $s_Y(x) \in \pi_{n+1}(Y)$
be the restriction to the second summand.
Clearly $s$ is natural with respect to $S^1$-equivariant maps.

The operation $s$ represents $\Delta$ on spherical classes, i.e.
the following diagram commutes for any $n>1$. 

\[\xymatrix{
\pi_{n+2}( X)\ar[d]\ar[r]^{s_X} & \pi_{n+3}(X)\ar[d]  \\
H_n(\Omega^2 X;k)\ar[r]^\Delta  & H_{n+1}(\Omega^2 X;k)
}\]

The only subtlety is
to show that $s_X$ can be replaced by $s_{\Omega^2 X}$.
Consider the universal example $X = S^n \wedge S^1_+$.
Clearly
$\pi_n(X)\ot k=H_n(X;k)=$\\
 $k<e_{n}>$,
and $\pi_{n+1}(X) \ot k = H_{n+1}(X;k)=k<e_{n+1}>$, because
$char(k)=0$ and $n>2$.
Moreover
$s_X(e_n)=e_{n+1}$ because 
the composite
$$S^1 \times S^n \hookrightarrow 
   S^1 \times X \to X \twoheadrightarrow 
   S^1 \wedge S^n = S^{n+1}$$
 has degree 1,  and $s_{\Omega^2 X}(e_n)=e_{n+1}$ for
the same reason.

Recall now that spherical classes are dual to indecomposables
in the minimal model.
In the operadic homology the BV-operator $t$  
on indecomposables in
$m(X)/(m(X) \cdot m(X))$
is induced by
the BV-operator on $m(X)$ by definition.
By naturality it is sufficient to show that $s$ and $t$ coincide in the
universal case
$X= S^n \wedge S^1_+$.
But $t(e_n)=e_{n+1}$, because, for homological reasons,
the map $f$ induced by the $S^1$-action on
 $m(X) = \Lambda(x_n,x_{n+1},\dots )$ gives
$f(x_{n+1})=x_{n+1} + x_n \ot e_1$.
\end{proof}


\begin{thebibliography}{1}
 

\bibitem{BV} J.M. Boardman and R.M. Vogt, {\em Homotopy invariant algebraic
structures on topological spaces}, LNM 347, 1973.
\bibitem{cohen} F. Cohen, Artin's braid groups, classical homotopy
theory, and sundry other curiosities, in 
{\em Braids}, Contemp. Math. 78, 167-206 (1988).
\bibitem{533} F. R. Cohen, T. J. Lada, J. P. May, {\em The homology of
iterated loop spaces}, LNM 533, 1976.
\bibitem{FeHaTh} Y. Felix, S. Halperin and
J.-C. Thomas, {\em Rational homotopy theory}, 
Rational homotopy theory, Graduate Texts in Mathematics, 205.
Springer-Verlag, New York, 2001.
\bibitem{Fiedo} Z. Fiedorowicz, The symmetric bar construction, preprint.
\bibitem{Fiedo2} Z. Fiedorowicz, Constructions of $E_n$ operads, preprint
math.AT/9808089, 1999.
\bibitem{getzler} E. Getzler, Batalin-Vilkovisky algebras and two-dimensional
topological field theories, {\em Comm. Math. Phys.} 159 n.2 (1994), 265-285. 
\bibitem{GJ} E. Getzler and J.D.S. Jones, Operads, homotopy algebra,
and iterated integrals for double loop spaces, preprint hep-th/9403055 (1994).
\bibitem{GK} V. Ginzburg and M. Kapranov, Koszul duality for operads, 
{\em Duke Math. Journal}  76, n.1, 203-272 (1994).
\bibitem{J} D. L. Johnson, Homeomorphisms of a surface which act
trivially on homology, {\em Proc. Math. American Soc.} 75, n.1,
119-125 (1979). 
\bibitem{JS} A. Joyal and R. Street, Braided tensor
categories, \emph{Advances in Mathematics}  102, 20-78 (1993).
\bibitem{muriel} M. Livernet, Homotopie rationnelle des alg\`ebres sur une 
op\'erade, PhD. Thesis, Strasbourg, 1998. 
\bibitem{Markl} M. Markl, Distributive laws and Koszulness, 
Ann. Inst. Fourier 46, n.2, 307-323 (1996).
\bibitem{Martin} M. Markl, A compactification of the real
configuration space as an operadic completion, {\em Journal of
Algebra} 215 (1999), 185-204.
\bibitem{book} M.Markl, S. Shnider and J. Stasheff,
Operads in algebra, topology and physics, book in preparation.
\bibitem{Milgram} R. James Milgram, Iterated loop spaces, 
{\em Ann. of Math.} (2) 84 (1966), 386--403.
\bibitem{Mimura} M. Mimura, H. Toda, {\em Topology of Lie groups I and II},
Transl. AMS 91, Providence 1991.
\bibitem{May}J.P. May, {\em The geometry of iterated loop spaces}, LNM
271, 1972.
\bibitem{Ray}N. Ray, The loop group of a mapping cone,
{\em Quart. J. Math.} 24 (1973), 485-498.
\bibitem{thesis} P. Salvatore, Configuration operads, minimal models
and rational curves, D. Phil. Thesis, Oxford, 1998.
\bibitem{barcelona} P. Salvatore, Configuration spaces with summable labels,
Proceedings of BCAT98, to appear in Progress in Math.
\bibitem{SV} R. Schw\"anzl and R.M. Vogt, {\em Coherence in homotopy group
actions}, LNM 1217, 364-390.
\bibitem{Vogt} R. Schw\"anzl and R.M. Vogt, The categories
of $A_\infty$ and $E_\infty$-monoids and ring spaces as closed simplicial
and topological categories, {\em Arch. Math.} 56 (1991), 405-411.
\bibitem{Segal} G. Segal, Configuration-spaces and iterated loop-spaces,
{\em Invent. Math.} 21 (1973), 213--221. 
\bibitem{SS} D. Sullivan, A. A. Voronov, paper in progress.
\bibitem{T} U. Tillmann, Higher Genus Surface Operad Detects
Infinite Loop Spaces, {\em Math. Ann.} 317 (2000), 613-628.
\bibitem{Vog} R. Vogt, Cofibrant Operads and Universal $E_\infty$ 
Operads, preprint.  
\bibitem{transfer} N. Wahl, Ribbon braid operad, Transfer Thesis,
Oxford 1999.
\bibitem{Nthesis} N. Wahl, Oxford University D.Phil. Thesis, in
preparation. 

\end{thebibliography}
\end{document}